\documentclass[]{amsart}
\usepackage{amsmath, amsthm, amscd, amsfonts, amssymb, graphicx, color}
\usepackage[bookmarksnumbered, colorlinks, plainpages]{hyperref}

\theoremstyle{definition}

\theoremstyle{remark}

\numberwithin{equation}{section}

\begin{document}
\setcounter{page}{1}
\begin{center}
{\bf   ARENS REGULARITY OF  BILINEAR FORMS AND\\ UNITAL BANACH MODULE SPACES}
\end{center}

\title[]{}

\author[]{KAZEM HAGHNEJAD AZAR  AND ABDOLHAMID RIAZI}

\address{}

\dedicatory{}



\begin{abstract}
Assume that $A$, $B$ are Banach algebras and $m:A\times B\rightarrow B$, $m^\prime:A\times A\rightarrow B$ are  bounded bilinear mappings. We will study the relation between Arens regularities of  $m$, $m^\prime$ and the Banach algebras $A$, $B$. For Banach $A-bimodule$ $B$, we show that $B$ factors with respect to $A$ if and only if $B^{**}$ is an unital $A^{**}-module$, and we define locally topological center for  elements of $A^{**}$ and will show that when locally topological center of mixed unit of $A^{**}$ is $B^{**}$, then $B^*$  factors on  both sides with respect to $A$ if and only if  $B^{**}$ has a  unit as $A^{**}-module$.
\end{abstract} \maketitle

\section{\bf  Preliminaries and
Introduction }

\noindent Throughout  this paper, $A$ is  a Banach algebra and $A^*$,
$A^{**}$, respectively, are the first and second dual of $A$. Recall that  a left approximate identity $(=LAI)$ [resp. right
approximate identity $(=RAI)$]
in Banach algebra $A$ is a net $(e_{\alpha})_{{\alpha}\in I}$ in $A$ such that   $e_{\alpha}a\longrightarrow a$ [resp. $ae_{\alpha}\longrightarrow a$]. We
say that a  net $(e_{\alpha})_{{\alpha}\in I}\subseteq A$ is a
approximate identity $(=AI)$ for $A$ if it is $LAI$ and $RAI$ for $A$. If $(e_{\alpha})_{{\alpha}\in I}$ in $A$ is bounded and $AI$ for $A$, then we say that $(e_{\alpha})_{{\alpha}\in I}$ is a bounded approximate identity  ($=BAI$) for $A$.
 For $a\in A$
 and $a^\prime\in A^*$, we denote by $a^\prime a$
 and $a a^\prime$ respectively, the functionals on $A^*$ defined by $<a^\prime a,b>=<a^\prime,ab>=a^\prime(ab)$ and $<a a^\prime,b>=<a^\prime,ba>=a^\prime(ba)$ for all $b\in A$.
 The Banach algebra $A$ is embedded in its second dual via the identification
 $<a,a^\prime>$ - $<a^\prime,a>$ for every $a\in
A$ and $a^\prime\in
A^*$.
  We denote the set   $\{a^\prime a:~a\in A~ and ~a^\prime\in
  A^*\}$ and
  $\{a a^\prime:~a\in A ~and ~a^\prime\in A^*\}$ by $A^*A$ and $AA^*$, respectively, clearly these two sets are subsets of $A^*$.\\
   Let $A$ have a $BAI$. If the
equality $A^*A=A^*,~~(AA^*=A^*)$ holds, then we say that $A^*$
factors on the left (right). If both equalities $A^*A=AA^*=A^*$
hold, then we say
that $A^*$  factors on both sides.\\
 The extension of bilinear maps on normed space and the concept of regularity of bilinear maps were studied by [1, 2, 5, 7, 12]. We start by recalling these definitions as follows.\\
 Let $X,Y,Z$ be normed spaces and $m:X\times Y\rightarrow Z$ be a bounded bilinear mapping. Arens in [1] offers two natural extensions $m^{***}$ and $m^{t***t}$ of $m$ from $X^{**}\times Y^{**}$ into $Z^{**}$ as following\\
1. $m^*:Z^*\times X\rightarrow Y^*$,~~~~~given by~~~$<m^*(z^\prime,x),y>=<z^\prime, m(x,y)>$ ~where $x\in X$, $y\in Y$, $z^\prime\in Z^*$,\\
 2. $m^{**}:Y^{**}\times Z^{*}\rightarrow X^*$,~~given by $<m^{**}(y^{\prime\prime},z^\prime),x>=<y^{\prime\prime},m^*(z^\prime,x)>$ ~where $x\in X$, $y^{\prime\prime}\in Y^{**}$, $z^\prime\in Z^*$,\\
3. $m^{***}:X^{**}\times Y^{**}\rightarrow Z^{**}$,~ given by~ ~ ~$<m^{***}(x^{\prime\prime},y^{\prime\prime}),z^\prime>$ \\$=<x^{\prime\prime},m^{**}(y^{\prime\prime},z^\prime)>$ ~where ~$x^{\prime\prime}\in X^{**}$, $y^{\prime\prime}\in Y^{**}$, $z^\prime\in Z^*$.\\
The mapping $m^{***}$ is the unique extension of $m$ such that $x^{\prime\prime}\rightarrow m^{***}(x^{\prime\prime},y^{\prime\prime})$ from $X^{**}$ into $Z^{**}$ is $weak^*-weak^*$ continuous for every $y^{\prime\prime}\in Y^{**}$, but the mapping $y^{\prime\prime}\rightarrow m^{***}(x^{\prime\prime},y^{\prime\prime})$ is not in general $weak^*-weak^*$ continuous from $Y^{**}$ into $Z^{**}$ unless $x^{\prime\prime}\in X$. Hence the first topological center of $m$ may  be defined as following
$$Z_1(m)=\{x^{\prime\prime}\in X^{**}:~~y^{\prime\prime}\rightarrow m^{***}(x^{\prime\prime},y^{\prime\prime})~~is~~weak^*-weak^*~~continuous\}.$$
Let now $m^t:Y\times X\rightarrow Z$ be the transpose of $m$ defined by $m^t(y,x)=m(x,y)$ for every $x\in X$ and $y\in Y$. Then $m^t$ is a continuous bilinear map from $Y\times X$ to $Z$, and so it may be extended as above to $m^{t***}:Y^{**}\times X^{**}\rightarrow Z^{**}$.
 The mapping $m^{t***t}:X^{**}\times Y^{**}\rightarrow Z^{**}$ in general is not equal to $m^{***}$, see [1], if $m^{***}=m^{t***t}$, then $m$ is called Arens regular. The mapping $y^{\prime\prime}\rightarrow m^{t***t}(x^{\prime\prime},y^{\prime\prime})$ is $weak^*-weak^*$ continuous for every $y^{\prime\prime}\in Y^{**}$, but the mapping $x^{\prime\prime}\rightarrow m^{t***t}(x^{\prime\prime},y^{\prime\prime})$ from $X^{**}$ into $Z^{**}$ is not in general  $weak^*-weak^*$ continuous for every $y^{\prime\prime}\in Y^{**}$. So we define the second topological center of $m$ as
$$Z_2(m)=\{y^{\prime\prime}\in Y^{**}:~~x^{\prime\prime}\rightarrow m^{t***t}(x^{\prime\prime},y^{\prime\prime})~~is~~weak^*-weak^*~~continuous\}.$$
It is clear that $m$ is Arens regular if and only if $Z_1(m)=X^{**}$ or $Z_2(m)=Y^{**}$. Arens regularity of $m$ is equivalent to the following
$$\lim_i\lim_j<z^\prime,m(x_i,y_j)>=\lim_j\lim_i<z^\prime,m(x_i,y_j)>,$$
whenever both limits exist for all bounded sequences $(x_i)_i\subseteq X$ , $(y_i)_i\subseteq Y$ and $z^\prime\in Z^*$, see [13].\\
The mapping $m$ is left strongly Arens irregular if $Z_1(m)=X$ and $m$ is right strongly Arens irregular if $Z_2(m)=Y$.\\
Let now $B$ be a Banach $A-bimodule$, and let\\
$$\pi_\ell:~A\times B\rightarrow B~~~and~~~\pi_r:~B\times A\rightarrow B.$$
be the left and right module actions of $A$ on $B$, respectively. Then $B^{**}$ is a Banach $A^{**}-bimodule$ with module actions
$$\pi_\ell^{***}:~A^{**}\times B^{**}\rightarrow B^{**}~~~and~~~\pi_r^{***}:~B^{**}\times A^{**}\rightarrow B^{**}.$$
Similarly, $B^{**}$ is a Banach $A^{**}-bimodule$ with module actions\\
$$\pi_\ell^{t***t}:~A^{**}\times B^{**}\rightarrow B^{**}~~~and~~~\pi_r^{t***t}:~B^{**}\times A^{**}\rightarrow B^{**}.$$
We may therefore define the topological centers of the left and right module actions of $A$ on $B$ as follows:\\
$$Z_{B^{**}}(A^{**})=Z(\pi_\ell)=\{a^{\prime\prime}\in A^{**}:~the~map~~b^{\prime\prime}\rightarrow \pi_\ell^{***}(a^{\prime\prime}, b^{\prime\prime})~:~B^{**}\rightarrow B^{**}$$$$~is~~~weak^*-weak^*~continuous\}$$
$$Z_{B^{**}}^t(A^{**})=Z(\pi_r^t)=\{a^{\prime\prime}\in A^{**}:~the~map~~b^{\prime\prime}\rightarrow \pi_r^{t***}(a^{\prime\prime}, b^{\prime\prime})~:~B^{**}\rightarrow B^{**}$$$$~is~~~weak^*-weak^*~continuous\}$$
$$Z_{A^{**}}(B^{**})=Z(\pi_r)=\{b^{\prime\prime}\in B^{**}:~the~map~~a^{\prime\prime}\rightarrow \pi_r^{***}(b^{\prime\prime}, a^{\prime\prime})~:~A^{**}\rightarrow B^{**}$$$$~is~~~weak^*-weak^*~continuous\}$$
$$Z_{A^{**}}^t(B^{**})=Z(\pi_\ell^t)=\{b^{\prime\prime}\in B^{**}:~the~map~~a^{\prime\prime}\rightarrow \pi_\ell^{t***}(b^{\prime\prime}, a^{\prime\prime})~:~A^{**}\rightarrow B^{**}$$$$~is~~~weak^*-weak^*~continuous\}$$

\noindent We note also that if $B$ is a left(resp. right) Banach $A-module$ and $\pi_\ell:~A\times B\rightarrow B$~(resp. $\pi_r:~B\times A\rightarrow B$) is left (resp. right) module action of $A$ on $B$, then $B^*$ is a right (resp. left) Banach $A-module$. \\
We write $ab=\pi_\ell(a,b)$, $ba=\pi_r(b,a)$, $\pi_\ell(a_1a_2,b)=\pi_\ell(a_1,a_2b)$, \\ $\pi_r(b,a_1a_2)=\pi_r(ba_1,a_2)$,~
$\pi_\ell^*(a_1b^\prime, a_2)=\pi_\ell^*(b^\prime, a_2a_1)$,~
$\pi_r^*(b^\prime a, b)=\pi_r^*(b^\prime, ab)$,~ \\for all $a_1,a_2, a\in A$, $b\in B$ and  $b^\prime\in B^*$
when there is no confusion.\\
Regarding $A$ as a Banach $A-bimodule$, the operation $\pi:A\times A\rightarrow A$ extends to $\pi^{***}$ and $\pi^{t***t}$ defined on $A^{**}\times A^{**}$. These extensions are known, respectively, as the first(left) and the second (right) Arens products, and with each of them, the second dual space $A^{**}$ becomes a Banach algebra. In this situation, we shall also simplify our notations. So the first (left) Arens product of $a^{\prime\prime},b^{\prime\prime}\in A^{**}$ shall be simply indicated by $a^{\prime\prime}b^{\prime\prime}$ and defined by the three steps:
 $$<a^\prime a,b>=<a^\prime ,ab>,$$
  $$<a^{\prime\prime} a^\prime,a>=<a^{\prime\prime}, a^\prime a>,$$
  $$<a^{\prime\prime}b^{\prime\prime},a^\prime>=<a^{\prime\prime},b^{\prime\prime}a^\prime>.$$
 for every $a,b\in A$ and $a^\prime\in A^*$. Similarly, the second (right) Arens product of $a^{\prime\prime},b^{\prime\prime}\in A^{**}$ shall be  indicated by $a^{\prime\prime}ob^{\prime\prime}$ and defined by :
 $$<a oa^\prime ,b>=<a^\prime ,ba>,$$
  $$<a^\prime oa^{\prime\prime} ,a>=<a^{\prime\prime},a oa^\prime >,$$
  $$<a^{\prime\prime}ob^{\prime\prime},a^\prime>=<b^{\prime\prime},a^\prime ob^{\prime\prime}>.$$
  for all $a,b\in A$ and $a^\prime\in A^*$.\\
  The regularity of a normed algebra $A$ is defined to be the regularity of its algebra multiplication when considered as a bilinear mapping. Let $a^{\prime\prime}$ and $b^{\prime\prime}$ be elements of $A^{**}$, the second dual of $A$. By $Goldstine^,s$ Theorem [6, P.424-425], there are nets $(a_{\alpha})_{\alpha}$ and $(b_{\beta})_{\beta}$ in $A$ such that $a^{\prime\prime}=weak^*-\lim_{\alpha}a_{\alpha}$ ~and~  $b^{\prime\prime}=weak^*-\lim_{\beta}b_{\beta}$. So it is easy to see that for all $a^\prime\in A^*$,
$$\lim_{\alpha}\lim_{\beta}<a^\prime,\pi(a_{\alpha},b_{\beta})>=<a^{\prime\prime}b^{\prime\prime},a^\prime>$$ and
$$\lim_{\beta}\lim_{\alpha}<a^\prime,\pi(a_{\alpha},b_{\beta})>=<a^{\prime\prime}ob^{\prime\prime},a^\prime>,$$
where $a^{\prime\prime}b^{\prime\prime}$ and $a^{\prime\prime}ob^{\prime\prime}$ are the first and second Arens products of $A^{**}$, respectively, see [5, 11, 12].\\
  We find the usual first and second topological center of $A^{**}$, which are
  $$Z_{A^{**}}(A^{**})=Z(\pi)=\{a^{\prime\prime}\in A^{**}: b^{\prime\prime}\rightarrow a^{\prime\prime}b^{\prime\prime}~ is~weak^*-weak^*$$$$~continuous\},$$
   $$Z^t_{A^{**}}(A^{**})=Z(\pi^t)=\{a^{\prime\prime}\in A^{**}: a^{\prime\prime}\rightarrow a^{\prime\prime}ob^{\prime\prime}~ is~weak^*-weak^*$$$$~continuous\}.$$
 An element $e^{\prime\prime}$ of $A^{**}$ is said to be a mixed unit if $e^{\prime\prime}$ is a
right unit for the first Arens multiplication and a left unit for
the second Arens multiplication. That is, $e^{\prime\prime}$ is a mixed unit if
and only if,
for each $a^{\prime\prime}\in A^{**}$, $a^{\prime\prime}e^{\prime\prime}=e^{\prime\prime}o a^{\prime\prime}=a^{\prime\prime}$. By
[3, p.146], an element $e^{\prime\prime}$ of $A^{**}$  is  mixed
      unit if and only if it is a $weak^*$ cluster point of some BAI $(e_\alpha)_{\alpha \in I}$  in
      $A$.\\
A functional $a^\prime$ in $A^*$ is said to be $wap$ (weakly almost
 periodic) on $A$ if the mapping $a\rightarrow a^\prime a$ from $A$ into
 $A^{*}$ is weakly compact. Pym in [13] showed that  this definition to the equivalent following condition\\
 For any two net $(a_{\alpha})_{\alpha}$ and $(b_{\beta})_{\beta}$
 in $\{a\in A:~\parallel a\parallel\leq 1\}$, we have\\
$$\lim_{\alpha}\lim_{\beta}<a^\prime,a_{\alpha}b_{\beta}>=\lim_{\beta}\lim_{\alpha}<a^\prime,a_{\alpha}b_{\beta}>,$$
whenever both iterated limits exist. The collection of all $wap$
functionals on $A$ is denoted by $wap(A)$. Also we have
$a^{\prime}\in wap(A)$ if and only if $<a^{\prime\prime}b^{\prime\prime},a^\prime>=<a^{\prime\prime}ob^{\prime\prime},a^\prime>$ for every $a^{\prime\prime},~b^{\prime\prime} \in
A^{**}$. \\
In all of this article, for two normed spaces $A$ and $B$, $\mathbf{B}(A,B)$ is the set of bounded linear operators from $A$ into $B$.\\
In next section, we will study the relationships between the Arens regularity of some bilinear mappings and Banach algebras, that is, for bounded bilinear mappings $m:A\times B\rightarrow B$, $m^\prime:A\times A\rightarrow B$, if $m$ or $m^\prime$ are Arens regular (resp. irregular), then $A$ or $B$ are Arens regular (resp. irregular), and conversely. We have some  applications from this discussion in some of algebra  as $L^1(G)$, $M(G)$, $L^{\infty}(X)$ and $C(X)$ whenever $G$ is a locally compact group and  $X$ is a semigroup. As a conclusion, with some conditions,  we show that if $A$ or $B$ is not Arens regular, then $A\hat{\otimes}B$ is not Arens regular.
 In chapter 3, we will extend some problems from [11] into module actions with some new results.\\ The main results of this paper can be summarized as follows:\\
{\bf a)} Let $A$ , $B$ be Banach algebras and $B$ be a Banach  $A-bimodule$. Let $T\in \mathbf{B}(A,B)$ be continuous and $m$ be the bilinear mapping from $A\times B$ into $B$ such that for every $a\in A$ and $b\in B$ we have $m(a,b)=T(a)b$. Then we have the following assertions\\
1) If $B$ is Arens regular, then $m$ is Arens regular.\\
2) If $T$ is surjective, then we have \\
i) $B$ is Arens regular if and only if $m$ is Arens regular.\\
ii) If $m$ is left strongly Arens irregular, then $B$ is left strongly Arens irregular.\\
iii) If $T$ is injective, then $B$ is left strongly Arens irregular if and only if $m$ is left strongly Arens irregular.\\
{\bf b)} Let $A$ , $B$ be Banach algebras and $B$ be a Banach $A-bimodule$. Let $T\in \mathbf{B}(A,B)$ be a   homomorphism. If $T$ is weakly compact, then the bilinear mapping $m(a_1,a_2)=T(a_1a_2)$ from $A\times A$ into $B$ is Arens regular.\\
{\bf c)} Assume that $B$ is a right Banach $A-module$ and $A^{**}$ has a right unit as $e^{\prime\prime}$.  Then,
  $B$  factors on the right with respect to $A$ if and only if  $e^{\prime\prime}$ is a right unit $A^{**}-module$ for $B^{**}$.\\
 {\bf d)} Assume that $A$ is a Banach algebra  and $A^{**}$ has a mixed unit $e^{\prime\prime}$.  Then we have the following assertions.\\
i) Let $B$ be a left Banach $A-module$. Then,  $B^*$  factors on the left with respect to $A$ if and only if $B^{**}$ has a left unit $A^{**}-module$ as $e^{\prime\prime}$.\\
ii) Let $B$ be a right Banach $A-module$ and $Z_{e^{\prime\prime}}(\pi_r^t)=B^{**}$. Then,  $B^*$  factors on the right with respect to $A$ if and only if $B^{**}$ has a right unit $A^{**}-module$ as $e^{\prime\prime}$.\\
iii) Let $B$ be a Banach  $A-bimodule$ and $Z_{e^{\prime\prime}}(\pi_r^t)=B^{**}$. Then,  $B^*$  factors on  both sides with respect to $A$ if and only if  $B^{**}$ has a  unit $A^{**}-module$ as $e^{\prime\prime}$.\\
\begin{center}
\section{ \bf Arens regularity of some bilinear forms }
\end{center}
In this part, we introduce some bilinear mappings from $A\times B$ or $A\times A$ into $B$ and make some relations between the Arens regularity of these  bilinear mappings and $A$ or $B$ with some applications.\\\\
{\it{\bf Theorem 1-2.}} Let $A$ , $B$ be Banach algebras and $B$ be a Banach  $A-bimodule$. Let $T\in \mathbf{B}(A,B)$ be continuous and $m$ be the bilinear mapping from $A\times B$ into $B$ such that for every $a\in A$ and $b\in B$ we have $m(a,b)=T(a)b$. Then we have the following assertions\\
a) If $B$ is Arens regular, then $m$ is Arens regular.\\
b) If $T$ is surjective, then we have \\
i) $B$ is Arens regular if and only if $m$ is Arens regular.\\
ii) If $m$ is left strongly Arens irregular, then $B$ is left strongly Arens irregular.\\
iii) If $T$ is injective, then $B$ is left strongly Arens irregular if and only if $m$ is left strongly Arens irregular.
\begin{proof}
 a) By  definition of $m^{***}$, we have $m^{***}(a^{\prime\prime},b^{\prime\prime})=T^{**}(a^{\prime\prime})b^{\prime\prime}$ and $m^{t***t}(a^{\prime\prime},b^{\prime\prime})=T^{**}(a^{\prime\prime})ob^{\prime\prime}$
where $a^{\prime\prime}\in A^{**}$ and $b^{\prime\prime}\in B^{**}$. Since $Z_1(B^{**})=B^{**}$,  the mapping $b^{\prime\prime}\rightarrow T^{**}(a^{\prime\prime})b^{\prime\prime}=m^{***}(a^{\prime\prime},b^{\prime\prime})$ is $weak^*-weak^*$ continuous for all $a^{\prime\prime}\in A^{**}$, also since $Z_2(B^{**})=B^{**}$, the mapping $a^{\prime\prime}\rightarrow T^{**}(a^{\prime\prime})ob^{\prime\prime}=m^{t***t}(a^{\prime\prime},b^{\prime\prime})$ is $weak^*-weak^*$ continuous for all $b^{\prime\prime}\in A^{**}$. Hence $m$ is Arens regular.\\
b) i) Let $m$ be Arens regular. Then $Z_1(m)=A^{**}$ and $Z_2(m)=B^{**}$. Let $b_1^{\prime\prime},b_2^{\prime\prime}\in B^{**}$ and $(b_{\alpha}^{\prime\prime})_{\alpha}\in B^{**}$ such that $b_{\alpha}^{\prime\prime}\stackrel{w^*} {\rightarrow}b_2^{\prime\prime}$. Assume that $a^{\prime\prime}\in A^{**}$ such that $T^{**}(a^{\prime\prime})=b_1^{\prime\prime}$. Then we have $$b_1^{\prime\prime}b_2^{\prime\prime}=T^{**}(a^{\prime\prime})b_2^{\prime\prime}=
m^{***}(a^{\prime\prime},b_2^{\prime\prime})=weak^*-\lim_{\alpha}m^{***}(a^{\prime\prime},b_{\alpha}^{\prime\prime})$$$$=
weak^*-\lim_{\alpha}T^{**}(a^{\prime\prime})b_{\alpha}^{\prime\prime}
=weak^*-\lim_{\alpha}b_1^{\prime\prime},b_{\alpha}^{\prime\prime}.$$
Hence $Z_1(B^{**})=B^{**}$  consequently $B$ is Arens regular.\\
b) ii) Let $m$ be left strongly Arens irregular then $Z_1(m)=A$. For  $b_1^{\prime\prime}\in Z_1(B^{**})$  the mapping $b_2^{\prime\prime}\rightarrow b_1^{\prime\prime}b_2^{\prime\prime}$~ is~ $weak^*-weak^*$ continuous. Also since $T$ is surjective,  there exists $a^{\prime\prime}\in A^{**}$ such that $T^{**}(a^{\prime\prime})=b_1^{\prime\prime}$ and the mapping
$b_2^{\prime\prime}\rightarrow T^{**}(a^{\prime\prime})b_2^{\prime\prime}=m^{***}(a^{\prime\prime},b_2^{\prime\prime})$ is $weak^*-weak^*$ continuous. Hence $a^{\prime\prime}\in Z_1(m)=A$. Consequently we have $b_1^{\prime\prime}=T^{**}(a^{\prime\prime})\in B$. It follows that $Z_1(B^{**})=B$.\\
b) iii) Let $B$ be left strongly Arens irregular, so $Z_1(B^{**})=B$. For $a^{\prime\prime}\in Z_1(m)$ the mapping $b^{\prime\prime}\rightarrow m^{***}(a^{\prime\prime},b^{\prime\prime})$ is $weak^*-weak^*$ continuous consequently
$T^{**}(a^{\prime\prime})\in Z_1(B^{**})=B$. Since $T$ is bijective,  $a^{\prime\prime}\in A$. Hence we conclude $Z_1(m)=A$.
\end{proof}
In Theorem 1-2, if we replace the  left strongly Arens irregularity of $A$, $B$ and $m$ with right strongly Arens irregularity of them, then the results will be similar.\\\\
The following definition which introduced by Ulger [17] has an important role in showing some sufficient condition for the Arens regularity of  tensor product $A\hat{\otimes} B$ where $A$ and $B$ are Banach algebra.\\
We recall that a bilinear form $m:A\times B\rightarrow C$ is biregular if for any two pairs of sequence
$({a}_i)_i$ ,  $(\tilde{a}_j)_j$ in $A_1$ and $(b_i)_i$, $(\tilde {b}_j)_j$
 in $B_1$, we have
$$\lim_i\lim_jm(a_i\tilde{a}_j,b_i\tilde{b}_j)=\lim_j\lim_im(a_i\tilde{a}_j,b_i\tilde{b}_j)$$
provided that these limits exist.\\\\
{\it{\bf Corollary 2-2.}} Let $B$ be an unital Banach algebra and suppose that $A$ is subalgebra of $B$.  If $A$ is not Arens regular, then $A\hat{\otimes}B$ is not Arens regular.
\begin{proof} Let  $m:A\times B\rightarrow C$ be the bilinear form that  introduced in Theorem 1-2 where $T:A\rightarrow B$ is natural inclusion. Since $A$ is not Arens regular,  $m$ is not biregular. Consequently by [17, Theorem 3.4], $A\hat{\otimes}B$ is not Arens regular.
\end{proof}
\noindent{\it{\bf Example 3-2.}} Let $X=[0,1]$ be the unit interval and let $C(X)$ be the Banach algebra of all continuous bounded functions on $X$ with supremum norm and the convolution as multiplication defined by
$$f*g(x)=\int_0^x f(x-t)g(t)dt~~~where~~ 0\leq x\leq 1.$$
Let $T:C(X)\rightarrow L^{\infty}(X)$ be the natural inclusion and $m:C(X)\times L^{\infty}(X)\rightarrow L^{\infty}(X)$ be defined by $m(f,g)=f*g$ where $f\in C(X)$ and $g\in L^{\infty}(X)$. By [2], $L^{\infty}(X)$ is Arens regular and by Theorem 1-2, we conclude that $m$ is Arens regular.\\
Similarly since $c_0$ is Arens regular, see[1, 5], by using Theorem 1-2, we can show that the bounded bilinear mapping $(f,g)\rightarrow f*g$ from $\ell^1\times c_0$ into $c_0$ is Arens regular.\\\\
For a Banach algebra $A$, we recall that a bounded linear operator $T:A\rightarrow A$ is said to be a left (resp. right) multiplier
if, for all $a, b\in A$, $T(ab)=T(a)b$ (resp. $T(ab)=aT(b)$). We denote by $LM(A)$ (resp. $RM(A)$) the set of all left (resp. right) multipliers of $A$. The set $LM(A)$ (resp. $RM(A)$) is normed subalgebra of the algebra $L(A)$ of bounded linear operator on $A$.\\
Now, we define a new concept as follows which is an extension  of  Left [right] multiplier on a Banach algebra. We will show some relation between this concept and Arens regularity of   some bilinear mappings in  Theorem 6-2. \\\\
{\it{\bf Definition 4-2.}} Let  $B$ be a left Banach [resp. right] $A-module$ and  $T\in \mathbf{B}(A,B)$. Then $T$ is called extended left [resp. right] multiplier if ~$T(a_1 a_2)=\pi_r(T(a_1),a_2)\\ ~[resp.~T(a_1 a_2)=\pi_\ell(a_1,T(a_2))]$ ~for all $a_1,a_2\in A$.\\
We show by $LM(A,B)~[resp. ~ RM(A,B)]$ all of the Left [resp. right] multiplier extension from $A$ into $B$. \\\\
{\it{\bf Example 5-2.}} Let $a^\prime\in A^*$. Then the mapping $T_{a^\prime}:a\rightarrow {a^\prime}a ~[resp.~~ R_{a^\prime}~a\rightarrow a{a^\prime}]$ from $A$ into $A^*$ is left [right] multiplier, that is, $T_{a^\prime}\in LM(A,A^*)~[R_{a^\prime}\in RM(A,A^*)].$ $T_{a^\prime}$ is weakly compact if and only if ${a^\prime}\in wap(A)$. So, we can write $wap(A)$ as a subspace of $LM(A,A^*)$.\\\\
{\it{\bf Theorem 6-2.}} Let $B$ be a left Banach $A-module$ and $T\in \mathbf{B}(A,B)$ be a continuous map. Assume that $m:A\times A\rightarrow B$ is the bilinear mapping such that $m(a_1,a_2)=T(a_1a_2)$. Then we have the following assertions\\
i) If $A$ is Arens regular, then $m$ is Arens regular.\\
ii) If $m$ is left [right] strongly Arens irregular, then $A$ is left [right] strongly Arens irregular.\\
iii) $T^{**}(Z_1(m))\subseteq Z_{A^{**}}(B^{**})$.\\
iv)  If $T\in LM(A,B)$, then $T^{**}\in LM(A^{**},B^{**})$ .\\
v) Suppose that $B$ is Banach algebra and $T$ is epimorphism. Then, $B$ is Arens regular if and only if $m$ is Arens regular.
\begin{proof} i) An easy calculation shows that $$m^{***}(a_1^{\prime\prime},a_2^{\prime\prime})=T^{**}(a_1^{\prime\prime}a_2^{\prime\prime})~~,~~ m^{t***t}(a_1^{\prime\prime},a_2^{\prime\prime})=T^{**}(a_1^{\prime\prime}oa_2^{\prime\prime}).$$ Since $A$ is Arens regular, the mapping $a_2^{\prime\prime}\rightarrow a_1^{\prime\prime} a_2^{\prime\prime}$ is $weak^*-weak^*$ continuous for all $a_1^{\prime\prime}\in A^{**}$. Also the mapping $a_1^{\prime\prime}\rightarrow a_1^{\prime\prime}o a_2^{\prime\prime}$ is $weak^*-weak^*$ continuous for all $a_2^{\prime\prime}\in A^{**}$. Hence  both mappings
$a_2^{\prime\prime}\rightarrow T^{**}(a_1^{\prime\prime}a_2^{\prime\prime})=m^{***}(a_1^{\prime\prime},a_2^{\prime\prime})$ and
$a_1^{\prime\prime}\rightarrow T^{**}(a_1^{\prime\prime}oa_2^{\prime\prime})=m^{t***t}(a_1^{\prime\prime},a_2^{\prime\prime})$ are $weak^*-weak^*$ continuous for all $a_1^{\prime\prime}\in A^{**}$ and $a_2^{\prime\prime}\in A^{**}$, respectively.
We conclude that $Z_1(m)=Z_2(m)=A^{**}$.\\
ii) Let $a_1^{\prime\prime}\in Z_1(A^{**})$. Then the mapping $a_2^{\prime\prime}\rightarrow a_1^{\prime\prime}a_2^{\prime\prime}$ is $weak^*-weak^*$ continuous consequently the mapping $a_2^{\prime\prime}\rightarrow T^{**}(a_1^{\prime\prime}a_2^{\prime\prime})=m^{***}(a_1^{\prime\prime},a_2^{\prime\prime})$ is $weak^*-weak^*$ continuous. Hence $a_1^{\prime\prime}\in Z_1(m)=A$.\\
iii) Let $a_1^{\prime\prime}\in Z_1(m)$. Then the mapping
$$a_2^{\prime\prime}\rightarrow m^{***}(a_1^{\prime\prime},
a_2^{\prime\prime})=T^{**}(a_1^{\prime\prime})a_2^{\prime\prime}$$ is $weak^*-weak^*$ continuous from $A^{**}$ into $B^{**}$. It follows that $T^{**}(a_1^{\prime\prime})\in Z_{A^{**}}(B^{**})$.\\
iv) If we set $m(a_1,a_2)=T(a_1a_2) [resp.~=T(a_1)a_2]$ for all $a_1,a_2\in A$, then $m^{***}(a_1^{\prime\prime},a_2^{\prime\prime})=T^{**}(a_1^{\prime\prime}a_2^{\prime\prime})~[resp. ~=T^{**}(a_1^{\prime\prime})a_2^{\prime\prime}]$ for all $a_1^{\prime\prime},a_2^{\prime\prime}\in A^{**}$. Thus, we conclude that $T^{**}(a_1^{\prime\prime}a_2^{\prime\prime})=T^{**}(a_1^{\prime\prime})a_2^{\prime\prime}$ for all $a_1^{\prime\prime},a_2^{\prime\prime}\in A^{**}$.\\
v) Let $m$ be Arens regular and $b_1^{\prime\prime},b_2^{\prime\prime}\in B^{**}$ and let $(b^{\prime\prime}_{\alpha})_{\alpha}\in B^{**}$ such that $b_{\alpha}^{\prime\prime}\stackrel{w^*} {\rightarrow}b_2^{\prime\prime}$. We set $a_1^{\prime\prime},a_2^{\prime\prime}\in A^{**}$ and $(a^{\prime\prime}_{\alpha})_{\alpha}\in A^{**}$ such that $T^{**}(a_1^{\prime\prime})=b_1^{\prime\prime}$ , $T^{**}(a_2^{\prime\prime})=b_2^{\prime\prime}$ and $T^{**}(a_{\alpha}^{\prime\prime})=b_{\alpha}^{\prime\prime}$. Then $$b_1^{\prime\prime}b_2^{\prime\prime}=T^{**}(a_1^{\prime\prime})T^{**}(a_2^{\prime\prime})
=T^{**}(a_1^{\prime\prime}a_2^{\prime\prime})
=m^{***}(a_1^{\prime\prime}
,a_2^{\prime\prime})$$
$$=weak^*-\lim_{\alpha}m^{***}(a_1^{\prime\prime}
,a_{\alpha}^{\prime\prime})=weak^*-\lim_{\alpha}T^{**}
(a_1^{\prime\prime}a_{\alpha}^{\prime\prime})$$
$$=weak^*-\lim_{\alpha}T^{**}
(a_1^{\prime\prime})T^{**}(a_{\alpha}^{\prime\prime})=
weak^*-\lim_{\alpha}b_1^{\prime\prime}b_{\alpha}^{\prime\prime},
$$
where by the open mapping Theorem, we have $a_{\alpha}^{\prime\prime}\stackrel{w^*} {\rightarrow}a_2^{\prime\prime}$. Consequently $Z_1(B^{**})=B^{**}$.\\
Conversely, let $B$ be Arens regular and $a_1^{\prime\prime},a_2^{\prime\prime}\in A^{**}$ and $(a_{\alpha})_{\alpha}\in A^{**}$ such that $a_{\alpha}^{\prime\prime}\stackrel{w^*} {\rightarrow}a_2^{\prime\prime}$.
 Then $$m^{***}(a_1^{\prime\prime}
,a_2^{\prime\prime})=T^{**}(a_1^{\prime\prime}a_2^{\prime\prime})=weak^*-\lim_{\alpha}T^{**}
(a_1^{\prime\prime}a_{\alpha}^{\prime\prime})$$$$=weak^*-\lim_{\alpha}m^{***}(a_1^{\prime\prime}
,a_{\alpha}^{\prime\prime}).$$ It follow that $Z_1(m)=A^{**}$. Thus $m$ is Arens regular.
\end{proof}
\noindent{\it{\bf Example 7-2.}} Assume that $T:c_0\rightarrow \ell^{\infty}$ is the natural inclusion and $m:c_0\times c_0\rightarrow \ell^{\infty}$ be the bilinear mapping such that $m(f,g)=f*g$. Since $c_0$ is Arens regular, then $m$ is Arens regular. Similarly the bilinear mapping $m:C(G)\times C(G)\rightarrow L^{\infty}(G)$ defined by formula $(f,g)\rightarrow f*g$ is Arens regular whenever $G$ is compact.\\\\
For normed spaces $X,Y,Z,W$, let $m_1:X\times Y\rightarrow Z$ and $m_2:X\times W\rightarrow Z$ be bounded bilinear mappings. If $h:Y\rightarrow W$ is a continuous linear mapping such that $m_1(x,y)=m_2(x,h(y))$ for all $x\in X$ and $y\in Y$, then we say that $m_1$ factors through $m_2$, see[2]. We say that the continuous bilinear mapping $m:X\times Y\rightarrow Z$  factors if $m$ is onto $Z$, see [7].\\\\
{\it{\bf Theorem 8-2.}} Let $A$ , $B$ be Banach algebras and $B$ be a Banach $A-bimodule$. Let $T\in \mathbf{B}(A,B)$ be a continuous  homomorphism. If $T$ is weakly compact, then the bilinear mapping $m(a_1,a_2)=T(a_1a_2)$ from $A\times A$ into $B$ is Arens regular.
\begin{proof} Let $m^\prime$ be the bilinear mapping that we introduced in Theorem 1-2. Then $m(a_1,a_2)=m^\prime(a_1,Ta_2)$ for all $a_1,a_2\in A$. Consequently $m$ factors through $m^\prime$, so by [2, Theorem 2], we conclude that $m$ is Arens regular.\\\\
\noindent{\it{\bf Example 9-2.}} Suppose that $T:L^1(G)\rightarrow M(G)$ is the natural inclusion. Then the bilinear mapping $m:L^1(G)\times L^1(G)\rightarrow M(G)$ defined by $m(f,g)=f*g$ for all $f,g\in L^1(G)$ is Arens regular whenever $G$ is finite, see [14]. Also the left strongly Arens irregularity of $m$ implies that $L^1(G)$ is also left strongly Arens irregular, see [10, 11].
\end{proof}
\section{ \bf Unital $A-modules$  and module actions}
\noindent In [11], Lau and Ulger show that for Banach algebra $A$, $A^*$ factors on the left if and only if $A^{**}$ is unital with respect to the first Arens product. In this chapter we extend this problem to module actions with some results.\\
We say that $A^{**}$ has a
$weak^*~bounded~left~approximate~identity (=W^*BLAI)$ with respect
to the first Arens product, if there is a bounded net as
$(e_{\alpha})_{\alpha}\subseteq A$ such that for all $a^{\prime\prime}\in
A^{**}$ and $a^{\prime}\in A^*$, we have
 $<e_{\alpha}a^{\prime\prime},a^{\prime}>\rightarrow<a^{\prime\prime},a^{\prime}>$. The definition of $W^*RBAI$ is similar to $W^*LBAI$   and
 if $A^{**}$ has both
$W^*LBAI$ and  $W^*RBAI$, then we say that $A^{**}$ has $W^*BAI$.\\
Assume that  $B$ is a Banach $A-bimodule$. We say that  $B$ factors on the left (right) with respect to $A$ if $B=BA~(B=AB)$. We say that $B$ factors on both sides, if $B=BA=AB$.\\\\
{\it{\bf Definition 1-3.}} Let $B$ be a left Banach $A-module$ and  $e$ be a left  unit element of $A$. Then we say that $e$ is a left unit (resp. weakly left unit)  $A-module$ for $B$, if $\pi_\ell(e,b)=b$ (resp. $<b^\prime , \pi_\ell(e,b)>=<b^\prime , b>$ for all $b^\prime\in B^*$) where $b\in B$. The definition of right unit (resp. weakly right unit) $A-module$ is similar.\\
We say that a Banach $A-bimodule$ $B$ is a unital $A-module$, if $B$ has left and right unit $A-module$ that are equal then we say that $B$ is unital $A-module$.\\
Let $B$ be a left Banach $A-module$ and $(e_{\alpha})_{\alpha}\subseteq A$ be a LAI [resp. weakly left approximate identity(=WLAI)] for $A$. We say that $(e_{\alpha})_{\alpha}$ is left approximate identity  ($=LAI$)[ resp. weakly left approximate identity  (=$WLAI$)] for $B$, if for all $b\in B$, we have $\pi_\ell (e_{\alpha},b) \stackrel{} {\rightarrow}
b$ ( resp. $\pi_\ell (e_{\alpha},b) \stackrel{w} {\rightarrow}
b$). The definition of the right approximate identity ($=RAI$)[ resp. weakly right approximate identity ($=WRAI$)] is similar.\\
We say that $(e_{\alpha})_{\alpha}$ is a approximate identity  ($=AI$)[ resp. weakly approximate identity  ($WAI$)] for $B$, if $B$ has left and right approximate identity  [ resp. weakly left and right approximate identity ] that are equal.\\
Let $(e_{\alpha})_{\alpha}\subseteq A$ be $weak^*$ left approximate identity for $A^{**}$. We say that $(e_{\alpha})_{\alpha}$ is $weak^*$ left approximate identity $A^{**}-module$ ($=W^*LAI~A^{**}-module$) for $B^{**}$, if for all $b^{\prime\prime}\in B^{**}$, we have $\pi_\ell^{***} (e_{\alpha},b^{\prime\prime}) \stackrel{w^*} {\rightarrow}
b^{\prime\prime}$. The definition of the $weak^*$ right approximate identity $A^{**}-module$($=W^*RAI~A^{**}-module$) is similar.\\
 We say that $(e_{\alpha})_{\alpha}$ is a $weak^*$ approximate identity $A^{**}-module$ ($=W^*AI~A^{**}-module$) for $B^{**}$, if $B^{**}$ has $weak^*$ left and right approximate identity $A^{**}-module$  that are equal.\\\\
{\it{\bf Example 2-3.}} i) $L^1(G)$ is a Banach $M(G)-bimodule$ under convolution as multiplication. It is clear that $L^1(G)$ is a $unital$  $M(G)-bimudole$. \\
ii) Since $L^p(G)$, for $1\leq p<\infty$, is a left Banach $M(G)-module$, by using [8, 10.15], $L^p(G)$ has a $BLAI$ $(e_{\alpha})_{\alpha}\subset M(G)$.\\\\
{\it{\bf Theorem 3-3.}} Assume that $A$ is a Banach algebra and has a $BAI$ $(e_{\alpha})_{\alpha}$. Then we have the following assertions.\\
i) Let $B$ be a right Banach $A-module$. Then $B$ factors on the left with respect to $A$ if and only if $B$ has a $WRAI$.\\
ii) Let $B$ be a left Banach $A-module$. Then $B$ factors on the right with respect to $A$ if and only if $B$ has a $WLAI$.\\
iii) $B$ factors on both side with respect to $A$ if and only if $B$ has a $WAI$.
\begin{proof} i) Suppose that $B=BA$. Let $b\in B$ and $b^\prime \in B^*$ then there are $x\in B$ and $a\in A$ such that $b=xa$. Then
$$<b^{\prime}, \pi_r (b,e_{\alpha})>=<b^{\prime}, \pi_r (xa,e_{\alpha})>=<\pi_r^*(b^{\prime}, x),ae_{\alpha}>\rightarrow<\pi_r^*(b^{\prime}, x),a>$$$$=<b^{\prime},\pi_r(x,a)>=<b^{\prime},b>.$$
It follows that $\pi_r (b,e_{\alpha}) \stackrel{w} {\rightarrow}
b$ consequently $B$ has a $WRAI$.\\
For the converse, since $BA$ is a weakly closed subspace of $B$, so by Cohen Factorization theorem, see [5], proof is hold.\\
ii) The proof is similar to (i).\\
iii) Clear.
\end{proof}
In Theorem 3-3, if we set $B=A$, then we obtain Lemma 2.1 from [11].\\\\
{\it{\bf Theorem 4-3.}}  Assume that $B$ is a right Banach $A-module$ and $A^{**}$ has a right unit as $e^{\prime\prime}$.  Then,
  $B$  factors on the right with respect to $A$ if and only if  $e^{\prime\prime}$ is a right unit $A^{**}-module$ for $B^{**}$.
\begin{proof} Since  $A^{**}$ has a right unit  $e^{\prime\prime}$, there is a $BRAI$ $(e_{\alpha})_{\alpha}$ for $A$ such that $e_{\alpha} \stackrel{w^*} {\rightarrow}e^{\prime\prime}$. Let $AB=B$ and $b\in B$. Thus, there is $x\in B$ and $a\in A$ such that $b=ax$. Then for all $b^\prime\in  B^*$, we have
$$<\pi_r^{**}(e^{\prime\prime},b^\prime),b>=<e^{\prime\prime},\pi_r^{*}(b^\prime,b)>
=\lim_{\alpha}<e_{\alpha},\pi_r^{*}(b^\prime,b)>$$
$$=\lim_{\alpha}<\pi_r^{*}(b^\prime,b),e_{\alpha}>=\lim_{\alpha}<b^\prime,\pi_r(b,e_{\alpha})>$$$$
=\lim_{\alpha}<b^\prime,\pi_r(xa,e_{\alpha})>=\lim_{\alpha}<\pi_r^{*}(b^\prime,x),ae_{\alpha}>$$
$$=<\pi_r^{*}(b^\prime,x),a>=<b^\prime, \pi_r(x,a)=<b^\prime,b>.$$
Thus $\pi_r^{**}(e^{\prime\prime},b^\prime)=b^\prime$. Now let $b^{\prime\prime}\in B^{**}$, then we have
$$<\pi_r^{***}(b^{\prime\prime},e^{\prime\prime}),b^\prime>
=<b^{\prime\prime},\pi_r^{**}(e^{\prime\prime},b^\prime)>=<b^{\prime\prime},b^\prime>.$$
We conclude that $\pi_r^{***}(b^{\prime\prime},e^{\prime\prime})=b^{\prime\prime}$.  Hence it follows that $B^{**}$ has a right unit $A^{**}-module$.\\
Conversely,   assume that  $e^{\prime\prime}$ is a right unit $A^{**}-module$ for $B^{**}$. Let $b\in B$ and $b^\prime\in B$. Then we have
$$<b^\prime,\pi_r(b,e_{\alpha})>=<\pi_r(b^\prime,b),e_{\alpha})>\rightarrow
<\pi_r(b^\prime,b),e^{\prime\prime}>=<b^\prime,\pi_r(b,e^{\prime\prime})>$$$$=<b^\prime,b>.$$
Consequently    $\pi_r(b,e_{\alpha})\stackrel{w} {\rightarrow}\pi_r(b,e^{\prime\prime})=b$,
it follows that $b\in \overline{BA}^{w}$. Since $BA$ is a $weakly$ closed subspace of $B$, so by Cohen Factorization theorem,  $b\in BA$.
\end{proof}
\noindent{\it{\bf Definition 5-3.}} Let $B$ be a Banach  $A-bimodule$ and $a^{\prime\prime}\in A^{**}$. We define the locally topological center of the left and right module actions of $a^{\prime\prime}$ on $B$, respectively, as follows\\
$$Z_{a^{\prime\prime}}^t(B^{**})=Z_{a^{\prime\prime}}^t(\pi_\ell^t)=\{b^{\prime\prime}\in B^{**}:~~~\pi^{t***t}_\ell(a^{\prime\prime},b^{\prime\prime})=
\pi^{***}_\ell(a^{\prime\prime},b^{\prime\prime})\},$$
$$Z_{a^{\prime\prime}}(B^{**})=Z_{a^{\prime\prime}}(\pi_r^t)=\{b^{\prime\prime}\in B^{**}:~~~\pi^{t***t}_r(b^{\prime\prime},a^{\prime\prime})=
\pi^{***}_r(b^{\prime\prime},a^{\prime\prime})\}.$$\\
It is clear that ~~~~~~~$$\bigcap_{a^{\prime\prime}\in A^{**}}Z_{a^{\prime\prime}}^t(B^{**})=Z_{A^{**}}^t(B^{**})=
Z(\pi_\ell^t),$$ ~~~~~~~ ~~\\ $$~~~~~~~~~~~~~~~~~~~~~~~~~~~~\bigcap_{a^{\prime\prime}\in A^{**}}Z_{a^{\prime\prime}}(B^{**})=Z_{A^{**}}(B^{**})=
Z(\pi_r).$$\\\\
{\it{\bf Theorem 6-3.}}  Assume that $A$ is a Banach algebra  and $A^{**}$ has a mixed unit $e^{\prime\prime}$.  Then we have the following assertions.\\
i) Let $B$ be a left Banach $A-module$. Then,  $B^*$  factors on the left with respect to $A$ if and only if $B^{**}$ has a left unit $A^{**}-module$ as $e^{\prime\prime}$.\\
ii) Let $B$ be a right Banach $A-module$ and $Z_{e^{\prime\prime}}(\pi_r^t)=B^{**}$. Then,  $B^*$  factors on the right with respect to $A$ if and only if $B^{**}$ has a right unit $A^{**}-module$ as $e^{\prime\prime}$.\\
iii) Let $B$ be a Banach  $A-bimodule$ and $Z_{e^{\prime\prime}}(\pi_r^t)=B^{**}$. Then,  $B^*$  factors on  both sides with respect to $A$ if and only if  $B^{**}$ has a  unit $A^{**}-module$ as $e^{\prime\prime}$.
\begin{proof} i) Let $(e_{\alpha})_{\alpha}\subseteq A$ be a BAI for $A$ such that  $e_{\alpha} \stackrel{w^*} {\rightarrow}e^{\prime\prime}$. Suppose that $B^*A=B^*$. Thus for all $b^{\prime}\in B^{*}$ there are $a\in A$ and $x^\prime\in B^*$ such that $x^\prime a=b^\prime$. Then for all $b^{\prime\prime}\in B^{**}$ we have
$$<\pi_\ell^{***}(e^{\prime\prime},b^{\prime\prime}),b^\prime>
=<e^{\prime\prime},\pi_\ell^{**}(b^{\prime\prime},b^\prime)>=
\lim_\alpha<\pi_\ell^{**}(b^{\prime\prime},b^\prime),e_{\alpha}>$$
$$=\lim_\alpha<b^{\prime\prime},\pi_\ell^{*}(b^\prime,e_{\alpha})>
=\lim_\alpha<b^{\prime\prime},\pi_\ell^{*}(x^\prime a,e_{\alpha})>$$
$$=\lim_\alpha<b^{\prime\prime},\pi_\ell^{*}(x^\prime ,ae_{\alpha})>
=\lim_\alpha<\pi_\ell^{**}(b^{\prime\prime},x^\prime) ,ae_{\alpha}>$$$$=<\pi_\ell^{**}(b^{\prime\prime},x^\prime) ,a>
=<\pi_\ell^{***}(b^{\prime\prime},b^{\prime}>.$$
Thus $\pi^{***}_\ell(e^{\prime\prime},b^{\prime\prime})=b^{\prime\prime}$ consequently $B^{**}$ has left unit $A^{**}-module$.\\
Conversely, Let $e^{\prime\prime}$ be a left unit $A^{**}-module$ for $B^{**}$ and $b^\prime\in B^*$. Then for all
$b^{\prime\prime}\in B^{**}$ we have
$$<b^{\prime\prime},b^\prime>=<\pi^{***}_\ell(e^{\prime\prime},b^{\prime\prime}),b^\prime>
 =<e^{\prime\prime},\pi^{**}_\ell(b^{\prime\prime},b^\prime)>$$
 $$=\lim_\alpha<\pi_\ell^{**}(b^{\prime\prime},b^\prime),e_{\alpha}>
 =\lim_\alpha<b^{\prime\prime},\pi_\ell^{*}(b^\prime,e_{\alpha})>.$$
 Thus we conclude that $weak-\lim_\alpha\pi_\ell^{*}(b^\prime,e_{\alpha})=b^\prime$. So by Cohen Factorization theorem, we are done.\\
ii) Suppose that $AB^*=B^*$. Thus for all $b^{\prime}\in B^{*}$ there are $a\in A$ and $x^\prime\in B^*$ such that $ax^\prime =b^\prime$. Assume $(e_{\alpha})_{\alpha}\subseteq A$ is a BAI for $A$ such that  $e_{\alpha} \stackrel{w^*} {\rightarrow}e^{\prime\prime}$. Let  $b^{\prime\prime}\in B^{**}$  and $(b_\beta)_\beta\subseteq B$ such that $b_\beta \stackrel{w^*} {\rightarrow}b^{\prime\prime}$. Then\\
$$<\pi_r^{***}(b^{\prime\prime},e^{\prime\prime}),b^\prime>=\lim_\beta<\pi_r^{***}(b_\beta,e^{\prime\prime}),b^\prime>=
 \lim_\beta \lim_\alpha<b^\prime,\pi_r(b_\beta,e_{\alpha})>$$
$$=\lim_\beta \lim_\alpha<ax^\prime,\pi_r(b_\beta,e_{\alpha})>=\lim_\beta \lim_\alpha<x^\prime,\pi_r(b_\beta,e_{\alpha})a>$$
$$=\lim_\beta \lim_\alpha<x^\prime,\pi_r(b_\beta,e_{\alpha}a)>=\lim_\beta \lim_\alpha<\pi_r^{*}(x^\prime,b_\beta),e_{\alpha}a)>$$
$$=\lim_\beta <\pi_r^{*}(x^\prime,b_\beta),a)>=<b^{\prime\prime},b^{\prime}>.$$
We conclude that $$\pi_r^{***}(b^{\prime\prime},e^{\prime\prime})=b^{\prime\prime}$$
for all $b^{\prime\prime}\in B^{**}$.\\
Conversely, Suppose that $\pi_r^{***}(b^{\prime\prime},e^{\prime\prime})=b^{\prime\prime}$ where $b^{\prime\prime}\in B^{**}$ and $(b_\beta)_\beta\subseteq B$ such that $b_\beta \stackrel{w^*} {\rightarrow}b^{\prime\prime}$. Let $(e_{\alpha})_{\alpha}\subseteq A$ be a BAI for $A$ such that  $e_{\alpha} \stackrel{w^*} {\rightarrow}e^{\prime\prime}$.  Then for all $b^\prime\in B^*$ we have
$$<b^{\prime\prime},b^{\prime}>=<\pi_r^{***}(b^{\prime\prime},e^{\prime\prime}),b^\prime>=
<b^{\prime\prime},\pi_r^{**}(e^{\prime\prime},b^\prime)>=\lim_\beta<\pi_r^{**}(e^{\prime\prime},b^\prime),b_\beta>$$
$$=\lim_\beta <e^{\prime\prime},\pi_r^{*}(b^\prime,b_\beta)>
=\lim_\beta \lim_\alpha<\pi_r^{*}(b^\prime,b_\beta),e_\alpha>$$
$$=\lim_\beta \lim_\alpha<\pi_r^{*}(b^\prime,b_\beta),e_\alpha>=\lim_\beta \lim_\alpha<b^\prime,\pi_r(b_\beta,e_\alpha)>$$
$$=\lim_\alpha \lim_\beta<\pi_r^{***}(b_\beta,e_\alpha),b^\prime>
=\lim_\alpha \lim_\beta<b_\beta,\pi_r^{**}(e_\alpha,b^\prime)>$$
$$= \lim_\alpha<b^{\prime\prime},\pi_r^{**}(e_\alpha,b^\prime)>.$$
It follows that $weak-\lim_\alpha\pi_r^{**}(e_\alpha,b^\prime)=b^\prime$. So by Cohen Factorization Theorem,  we are done.\\
iii) Clear.
\end{proof}
\noindent{\it{\bf Corollary 7-3.}} Let $B$ be a Banach  $A-bimodule$ and $A^{**}$ has a mixed unit $e^{\prime\prime}$.  \\
a) Let $Z_{e^{\prime\prime}}(\pi_r^t)=B^{**}$. Then we have the following assertions\\
i) If $B$ or $B^*$ factors on the right but not on the left with respect to $A$ then $\pi_\ell\neq\pi_r^t$.\\
ii) If  $B^*$ factors on the left  with respect to $A$ and $\pi_\ell=\pi_r^t$ , then  $B^*$ factors on the right with respect to $A$.\\
iii) $e^{\prime\prime}$ is a right unit $A^{**}-module$ for $B^{**}$ if and only if $(e_{\alpha})_{\alpha}$ is a $W^*RAI$ $A^{**}-module$ for $B^{**}$ whenever $e_{\alpha} \stackrel{w^*} {\rightarrow}e^{\prime\prime}$.\\
b) Let $Z_{e^{\prime\prime}}^t(\pi_\ell^t)
=B^{**}$. Then we have the following assertions\\
i) If $B$ or $B^*$ factors on the right but not on the left with respect to $A$ then $\pi_r\neq\pi_\ell^t$.\\
ii) If  $B^*$ factors on the right  with respect to $A$ and $\pi_r=\pi_\ell^t$ , then  $B^*$ factors on the left with respect to $A$.\\
iii) $e^{\prime\prime}$ is a left unit $A^{**}-module$ for $B^{**}$ if and only if $(e_{\alpha})_{\alpha}$ is a $W^*LAI$ $A^{**}-module$ for $B^{**}$ whenever $e_{\alpha} \stackrel{w^*} {\rightarrow}e^{\prime\prime}$.\\
c) Let $Z_{e^{\prime\prime}}^t(\pi_\ell^t)
= Z_{e^{\prime\prime}}(\pi_r^t)=B^{**}$. Then we have the following assertions\\
i) If  $B^*$ not factors on the right and left with respect to $A$ then $\pi_r\neq\pi_\ell^t$ and $\pi_\ell\neq\pi_r^t$.\\
ii) $e^{\prime\prime}$ is a unit $A^{**}-module$  for $B^{**}$ if and only if $(e_{\alpha})_{\alpha}$ is a $W^*AI$ $A^{**}-module$ for $B^{**}$ whenever $e_{\alpha} \stackrel{w^*} {\rightarrow}e^{\prime\prime}$.
\begin{proof} a) i)  Let $B$ or $B^*$ factors on the right but not on the left with respect to $A$. By Theorem 4-3 (resp. Theorem 5-3),  $e^{\prime\prime}$ is a right unit $A^{**}-module$ for $B^{**}$. Thus we have $\pi_r^{***}(b^{\prime\prime},e^{\prime\prime})=b^{\prime\prime}$ for all $b^{\prime\prime}\in B^{**}$. If we set $\pi_\ell=\pi_r^t$, then $\pi_\ell^{***}(e^{\prime\prime},b^{\prime\prime})=\pi_r^{t***}(e^{\prime\prime},b^{\prime\prime})
=\pi_r^{t***t}(b^{\prime\prime},e^{\prime\prime})=\pi_r^{***}(b^{\prime\prime},e^{\prime\prime})=b^{\prime\prime}$ for all $b^{\prime\prime}\in B^{**}$. Consequently $e^{\prime\prime}$ is left unit $A^{**}-module$ for $B^{**}$. Then by Theorem 4-3 (resp. Theorem 5-3), $B$ or $B^*$ factors on the left which is impossible.\\
 ii) Similar to (i).\\
 iii) Since $e_{\alpha} \stackrel{w^*} {\rightarrow}e^{\prime\prime}$, $weak^*-\lim_{\alpha}\pi _r^{***}(b^{\prime\prime},e_{\alpha})=\pi _r^{***}(b^{\prime\prime},e^{\prime\prime})$ for all $b^{\prime\prime}\in B^{**}$ hence we are done.\\
The proofs of (b) and (c) is the same and easy.
\end{proof}
 Assume that  $Z_{e^{\prime\prime}}^t(\pi_\ell^t)
= Z_{e^{\prime\prime}}(\pi_r^t)=B^{**}$. Let $\pi_r=\pi_\ell^t$ and $\pi_\ell=\pi_r^t$. We  conclude from Corollary 7-3 that if  $B^*$ also factors on the one side, then  $B^*$ factors on the other side.\\\\
 In Theorem 7-3. if we set $B=A$, then we obtain the Proposition 2.10 from [11].\\\\
 {\it{\bf Problems.}} \\
 Suppose that $B$ is a  Banach $A-bimodule$. Which condition we need for the following assertions\\
 i) ~~ $B$  factors on the left with respect to $A$ if and only if $B^{**}$ has a left unit $A^{**}-module$.\\
 ii)  $B$  factors on the one side with respect to $A$ if and only if  $B^*$  factors on the same side with respect to $A$.

\bibliographystyle{amsplain}

\it{Department of Mathematics, Amirkabir University of Technology, Tehran, Iran\\
{\it Email address:} haghnejad@aut.ac.ir\\\\
Department of Mathematics, Amirkabir University of Technology, Tehran, Iran\\
{\it Email address:} riazi@aut.ac.ir}
\end{document}